\documentclass[11pt]{article} 
\usepackage{amsfonts}
\usepackage[utf8]{inputenc} 
\usepackage{geometry} 
\usepackage{graphicx} 

\usepackage{booktabs} 
\usepackage{array} 
\usepackage{paralist} 
\usepackage{verbatim} 
\usepackage{subfig} 
\usepackage{biblatex}
\usepackage{csquotes}
\usepackage[english]{babel}
\usepackage{amsmath}
\usepackage[utf8]{inputenc}
\usepackage{graphicx} 
\usepackage{booktabs} 
\usepackage{mathrsfs}
\usepackage{fancyhdr} 
\usepackage{amsmath}
\usepackage{sectsty}
\allsectionsfont{\sffamily\mdseries\upshape} 

\usepackage[nottoc,notlof,notlot]{tocbibind} 
\usepackage[titles,subfigure]{tocloft} 
\usepackage{indentfirst}

\usepackage{amsthm}

\theoremstyle{plain}
\newtheorem{theorem}{Theorem}[subsection]
\newtheorem{lemma}{Lemma}[subsection]
\newtheorem{proposition}{Proposition}[subsection]
\newtheorem{corollary}{Corollary}[subsection]
\theoremstyle{definition}
\newtheorem{definition}{Definition}[subsection]
\newtheorem{conjecture}{Conjecture}
\newtheorem{example}{Example}[subsection]
\theoremstyle{remark}

\title{On the Analysis of Boolean Functions and Fourier-Entropy-Influence Conjecture}
\author{Xiao Han}
\begin{document}
\date{}
\thispagestyle{empty}
\maketitle
\thispagestyle{empty}
\newpage
\setcounter{page}{1}
\tableofcontents
\newpage
\section{Introduction}
\subsection{About boolean functions}
\indent Analysis of Boolean Functions is part of a lot of subjects including (Discrete) Harmonic Analysis,  Combinatorics, Probability and Theoretical Computer Science due to numerous applications of it in all these fields. For example, it has strong connections with several fashionable topics like learning theory, percolation theory and random graph. The definition of boolean function is simple:
$$
f:\{-1,1\}^n \to \{-1,1\} \text{ is called a boolean function}.
$$
\indent Despite of the simple definition, boolean function actually exhibits large amount of untrivial characteristics. The method of Fourier-Walsh expansion turns out to be indispensable in the study of boolean functions. Coming from pure functional analysis, it also works well in the discrete setting  and becomes one of the most important components in the analysis of boolean functions.\\
\indent A core topic in the analysis of boolean functions is the study of boolean functions with low influence. Although we already have a lot of powerful results like KKL theorem and Friedgut's junta theorem to characterize them, it seems that we still know little in this topic. The famous Fourier-Entropy-Influence conjecture is a typical example. The slow progress of it despite of significant efforts seems to imply the lack of powerful tools.\\
\subsection{About Fourier-Entropy-Influence conjecture}
\indent The Fourier-Entropy-Influence conjecture asks if the spectrum entropy of a boolean function can be bounded by a constant times of its influence. Or more precisely, if there exists a constant $c>0$ such that
$$
Ent(f)\leq c I(f) ?
$$
\indent One may see Section 2.4 for a more rigorous statement. In Section 2.4 we also show that the FEI conjecture implies that the Fourier weights of boolean functions are highly concentrated on $\exp(O(I(f)))$ coefficients, which will be useful in learning theory, random graph and other fields. The famous Friedgut's junta theorem which describes the structure of Fourier spectrum of boolean functions could only give a result of $\exp(O(I(f)^2))$ (see Section 3.4 for more infomation), thus we could see that efforts towards the FEI conjecture may lead us to understand more about the structrue of Fourier spectrum of boolean functions.\\
\indent The conjecture is first given by Friedgut and Kalai in order to study the threshold phenomena of monotone graph properties in random graphs  \cite{10}. As we said, the progress towards the FEI conjecture has in fact been slow and it is only proved for certain classes of boolean functions. (See \cite{1}, \cite{8}, \cite{12}, \cite{13} for example.) In a recent paper they give a slightly weaker bound for Fourier entropy of the low-degree part of boolean functions \cite{6}.\\
\indent A strictly weaker version of the FEI conjecture is the Fourier-Min-Entropy-Influence conjecture, which asks if the maximum value of the Fourier coefficients of a boolean function $f$ can be lower bounded by $\exp(-O(I(f)))$. (For a more formal statement also see Section 2.4.) Even this conjecture is hard to be resolved. By Friedgut's junta theorem we can get a bound of $\exp(-O(I(f)^2))$ and in \cite{6} they gave a better bound of $\exp(-O(\widetilde{I}(f)\log(\widetilde{I}(f)+1)))$ where $\widetilde{I}(f)=\frac{I(f)}{Var(f)}$ (or $\exp(-O(I(f)\log(I(f)+1)))$ approximately).\\
\subsection{About this work}
\indent As a master thesis, this manuscript also includes some classical and meaningful results we select apart from the new results we've found. We try to ensure the self-completeness of this work so that readers could probably read it independently. The purpose of this work is to summarize part of the untrivial existing results in the fields of Analysis of Boolean Functions and try to develop them by some new theories and applications that we've found.\\
\indent The results in this work which are found by us are Theorem 4.2.3, Theorem 4.3.1 and Theorem 5.2.1-5.2.3. Additionally we give Proposition 4.4.1 and Proposition 4.4.2 after Theorem 4.3.1 to get some examples about it. Perhaps among them the most meaningful one is Theorem 5.2.2, which is stated as follows: \\
\textbf{Theorem 5.2.2.} \emph{There exist $c_1, c_2>0$, such that for any boolean function $f:\{-1,1\}^n\to \{-1,1\}$ we have}
$$Ent(f)\leq c_1 I(f)+ c_2 \sum_{k\in[n]}-I_k(f)\log I_k(f).$$
Theorem 5.2.2 is a rather good new estimation for $Ent(f)$. In the proof we actually find some bounds for $\sum_S \hat{f}(S)^{2(1+\epsilon)}$ by the influence, where $\epsilon\in (0,\frac{1}{2})$. This could probably be improved as well and also somehow implies the extent of the concentration of Fourier coefficients (or even the distribution of the values of Fourier coefficients), which might help us improve our understanding on the structure of Fourier spectrum.\\
\indent We introduce the construction of this work as the end of the introduction. In Chapter 2 we give the basic definitions about boolean functions, with a more concrete example in Section 2.3 to show its application and introduce the FEI conjecture and FMEI conjecture. In Chapter 3 we concentrate on the classical topic in analysis of boolean functions: the hypercontractivity and its application. As for Chapter 4, we'll show another object, High-order Influence, which is long known but there's not much results about it in the literature. For it we give some existing results, show our own development of it and some applications that we've found as well. At last in Chapter 5 we introduce the idea of restrictions in dealing with boolean functions and use it to prove Theorem 5.2.3.

\section{Preliminaries}
\subsection{Definition of boolean function and Fourier-Walsh expansion}
\indent For $n\in \mathbb{N}^*$, let $\Omega_n^*$ denotes the function space of all the function from $\{-1,1\}^n$ to $\mathbb{R}$, boolean functions are in fact, a certain class of functions in $\Omega_n^*$. We first give the definition.
\begin{definition}
A function $f: \{-1,1\}^n \to \{-1,1\}$ from the vertex of an hypercube to a two-point set is called a boolean function.
\end{definition}
Endow $\{-1,1\}^n$ with a uniform probability measure $\mathbb{P}$, so that $x\in \{-1,1\}^n$ becomes a uniform random vector in this probability space and we can therefore define the $L^p$ norms on $\Omega_n^*$: for any $f\in \Omega_n^*$, $p\geq 1$, $||f||_p:= (\mathbb{E}_x f(x)^p)^{\frac{1}{p}}$. We also get a natural inner product on $\Omega_n^*$: for any $f,g\in \Omega_n^*$,
\begin{equation}
<f,g>:=\mathbb{E}_x f(x)g(x).
\end{equation}
\indent Furthermore, the classical functional analysis helps us give the Fourier-Walsh expansion of functions in $\Omega_n^*$. For any $x=(x_1,x_2,...,x_n)\in \{-1,1\}^n$, if we define $X_S(x):=\prod_{i\in S}x_i$, where $S$ is any subset of $[n]=\{1,2,...,n\}$ and in particular $X_\emptyset=1$, we'll have
\begin{theorem}
$ \{X_S\}$ is an orthogonal basis of $\Omega_n^*$, so that for any $f:\{-1,1\}^n \to \mathbb{R}$, we can write $f(x)=\sum_{S\in[n]}{\hat{f}(S)X_S}$(the Fourier-Walsh expansion), where $\hat{f}(S)=\frac{1}{2^n}\sum_{x\in \{-1,1\}^n}f(x)X_S(x)$ are called the Fourier coefficients of $f$.
\end{theorem}
\textbf{Proof} For orthogonality, we only need to check that for $S \neq T$, $\sum_{x}X_S X_T=0$. To show this we simply need to find that for $i \in S$ but $i \notin T$, $X_S X_T$ changes its signal when $x_i$ changes its signal. Note that $X_S$ are nonzero functions, by calculating the dimension we've already show that $\{X_S\}$ is an orthogonal basis of $V$. On the other hand, by orthogonality we have $\sum_{x}f(x)X_S= \sum_{x} \sum_{T\in[n]} \hat{f}(S) X_S X_T= \sum_{x} \hat{f}(S) X_S^2= 2^n \hat{f}(S)$, which helps us get the expressions of Fourier coefficients. \qed.\\\\
\indent \textbf{Remark} By the Fourier-Walsh expansion, we actually write $f$ into the form of a polynomial of $x$ where $\hat{f}(S)$ are the coefficients of this polynomial.\\ \\
\indent One may also check that $||X_S||_2^2=\mathbb{E}X_S^2=1$ to show the basis is also normalized. Thus we have Parseval's formula:
\begin{theorem}
For any $f: \{-1,1\}^n \to \mathbb{R}$, $||f||_2^2=\mathbb{E}f(x)^2=\sum_{S \in [n]}\hat{f}(S)^2$.
\end{theorem}
Note that for boolean function $f$, $f(x)^2=1$ for any $x$, so a direct calculation shows that $\sum_{S \in [n]}\hat{f}(S)^2=\mathbb{E}f(x)^2=1$. From now on, if there's no other explanation, we always asume that $x$ obeys such uniform distribution when it is considered as a random vector.\\
\indent Under such assumption, $\mathbb{E}f=\hat{f}(\emptyset)$ so that
\begin{equation}
Var (f)=\mathbb{E}f^2-(\mathbb{E}f)^2=1- \hat{f}(\emptyset)^2=\sum_{|S|\geq1}\hat{f}(S)^2.
\end{equation}
\indent At last, we give the definition of the spectral sample:
\begin{definition}
For a boolean function $f$, its spectral sample $\mathcal{S}_f$ is a random variable takes value in $\mathscr{F}_n=\{S : S \subseteq [n]\}$ such that $\mathbb{P}(\mathcal{S}_f=S)=\hat{f}(S)^2$ for any $S\in \mathscr{F}_n$.
\end{definition}
\subsection{Pivotality and influence}
\indent For any $x=(x_1,x_2,...x_n) \in \{-1,1\}^n, k\in [n]$, we write $\mu_k(x)=(x_1,..,x_{k-1},-x_k,x_{k+1},..,x_n)$ to denote the point generated by flipping the $k^{th}$ coordinate of $x$.\\
\indent We first give the definition of pivotality.
\begin{definition}
For a boolean function $f$, $k\in [n]$ is pivotal for $x\in \{-1,1\}^n$ if and only if $f(x)\neq f(\mu_k(x))$.
\end{definition}
\indent From the definition we may see that $k$ is pivotal for $x$ basically means that $x_k$ has some influences on the quantity of $f(x)$. We furthermore give the definition of pivotal set.
\begin{definition}
For a boolean function $f$, $\mathcal{P}(x):=\{k \in [n] : k \text{ is pivotal for }x\}$ is called the pivotal set of $x$ for $f$.
\end{definition}
\indent In short, pivotal set consists of all those bits which are pivotal. We also directly give the definitions of influence and total influence.
\begin{definition}
For a boolean function $f$, the influence of the $k^{th}$ bit is defined by $I_k(f):=\mathbb{P}(k \in \mathcal{P}(x))$.
\end{definition}
\begin{definition}
For a boolean function $f$, the total influence of $f$ is defined by $I(f):=\sum_{k\in[n]}I_k(f)$.
\end{definition}
Note that we have $I(f)= \mathbb{E}_x |\mathcal{P}(x)|$ directly from the definitions. The theorem below shows the connection between the influence and the Fourier coefficients of $f$.
\begin{theorem}
Let $f$ be a boolean function, $I_k(f)=\sum_{S \subseteq [n], k\in S}\hat{f}(S)^2$ for any $k\in [n]$ and $I(f)=\sum_{S \subseteq [n]}|S|\hat{f}(S)^2$.
\end{theorem}
\textbf{Proof} Note that $I_k(f)= \mathbb{E}_x|\frac{1}{2}(f(x)-f(\mu_k(x)))|=\mathbb{E}_x(\frac{1}{2}(f(x)-f(\mu_k(x))))^2$. Under Fourier-Walsh expansion, we have $f(x)=\sum_{S\in [n]}\hat{f}(S)X_S$ and $f(\mu_k(x))=\sum_{k\notin S}{\hat{f}(S)X_S}-\sum_{k\in S}{\hat{f}(S)X_S}$. So $I_k(f)=\mathbb{E}_x(\sum_{k\in S}\hat{f}(S)X_S)^2$. By Parseval's formula, we have $I_k(f)=\sum_{k\in S}\hat{f}(S)^2$ and $I(f)=\sum_{k\in[n]}I_k(f)=\sum_{S \subseteq [n]}|S|\hat{f}(S)^2$. \qed\\\\
\indent By this theorem we have $I(f)=\mathbb{E}|\mathcal{S}_f|$ for any boolean function $f$.
\subsection{General situation for independent bits and percolation}
We now focus on a concrete random model. Let $y=(y_1,y_2,...,y_n) \in \{-1,1\}^n$, where $y_1, y_2,...y_n$ are independent random variables such that $\mathbb{P}(y_i=1)=\frac{1+p_i}{2}$, where $p_i \in[-1,1]$. $f$ is a boolean function. For convenience we write $p_S:=\prod_{i\in S}p_i$. Note that $\mathbb{P}(y_i=-1)=\frac{1-p_i}{2}$, so that we have for any $x=(x_1,x_2,...,x_n)\in \{-1,1\}^n$
\begin{equation}
\mathbb{P}(y_i=x_i)=\frac{1+x_i p_i}{2},
\end{equation}
thus,
\begin{equation}
\begin{aligned}
\mathbb{E}_y f(y)&=\sum_{x\in \{-1,1\}^n}\mathbb{P}(y=x)f(x)\\
&=\sum_{x\in \{-1,1\}^n}(2^{-n}f(x)\prod_{i=1}^{n}(1+x_i p_i))\\
&=\sum_{x\in \{-1,1\}^n}{(2^{-n}f(x)\sum_{S\subseteq [n]}p_S X_S)}\\
&=\sum_{S\subseteq [n]}(2^{-n}\sum_{x\in \{-1,1\}^n}f(x)X_S)p_S\\
&=\sum_{S\subseteq [n]}\hat{f}(S)p_S.
\end{aligned}
\end{equation}
\indent Note that $\mathbb{E}_y f(y)=2\mathbb{P}(f(y)=1)-1$, (2.4) can in fact be seen as a Taylor expansion of $\mathbb{P}(f(y)=1)$ at the point $\mathbb{P}(y_i=1)=\frac{1}{2}$ for any $i\in [n]$.
\begin{example}
\textbf{(Percolation)} If $y$ is a sample point in the percolation model, so that $\mathbb{P}(y_i=1)=p$ for any $i\in [n]$. We define $f_A(y):= 1_{y\in A}-1_{y \notin A}$ for any event $A$ where $f_A$ is a boolean function, then from (2.4) we'll have
\begin{equation}
\frac{d(\mathbb{P}(A))}{dp}(\frac{1}{2})=\sum_{i\in[n]}\widehat{f_A}(\{i\}).
\end{equation}
\indent Notice that for increasing events A, $f_A$ is a monotone function and $I_k(f_A)=\mathbb{E}_x |\frac{1}{2}(f_A(x)-f_A(\mu_k(x)))|=\mathbb{E}_x f_A(x)x_k=\widehat{f_A}(\{k\})$, thus $\frac{d(\mathbb{P}(A))}{dp}(\frac{1}{2})=\sum_{k\in[n]}I_k(f_A)=I(f_A)$, which is the famous Margulis-Russo formula at $p=\frac{1}{2}$.
\end{example}
\subsection{Fourier-Entropy-Influence conjecture}
\indent The Fourier-Entropy-Influence Conjecture (or FEI conjecture in short) asks if the entropy of the Fourier coefficients of a boolean function can always be bounded by a constant times of its influence. Here the entropy is defined as follows.
\begin{definition}
For any boolean function $f$, the spectrum entropy of it is defined as $Ent(f):=\sum_{S}\hat{f}(S)^2\log_2 \frac{1}{\hat{f}(S)^2}$.
\end{definition}
Note that it is slightly different from the classical definition. And the FEI conjecture is:
\begin{conjecture}
(\textbf{FEI})  There exists a constant $c>0$, such that for any boolean function $f$, we have $Ent(f) \leq c I(f)$.
\end{conjecture}
\indent Note that if FEI conjecture is right, we have for any $\delta>0$, $\sum_{\hat{f}(S)^2\leq 2^{-\frac{cI(f)}{\delta}}}\hat{f}(S)^2\leq \delta$, on the other hand $|\{S:\hat{f}(S)^2> 2^{-\frac{cI(f)}{\delta}}\}|\leq 2^{\frac{cI(f)}{\delta}}$, which means that the Fourier weights are actually concentrated on $\exp(O(I(f)))$ coefficients except a constant of $\delta$.\\
\indent The following example shows that such property doesn't always hold for non-boolean functions.
\begin{example}
For $f:\{-1,1\}^n \to \mathbb{R}$, let $f=\sum_S \hat{f}(S)X_S$, where $\hat{f}(S)=1_{|S|=1}\cdot n^{-\frac{1}{2}}$. One may check that $\sum_S |S|\hat{f}(S)^2=1$ (although we didn't define influence for non-boolean functions, we use the fact that $I(f)=\sum_S |S|\hat{f}(S)^2$ here instead). On the other hand, when $n\to \infty$, $Ent(f)=\log_2 n \to \infty$, thus we couldn't find the constant $c$ for such $f$.
\end{example}
\indent Now we introduce the FMEI conjecture. If the FEI conjecture is right, since $\sum_{S}\hat{f}(S)^2=1$ we have $\min \log_2 \frac{1}{\hat{f}(S)^2} \leq c I(f)$, which is equivalent to the following Fourier-Min-Entropy-Influence Conjecture.
\begin{conjecture}
(\textbf{FMEI})  There exists a constant $c>0$, such that for any boolean function $f$, we have $\max_S |\hat{f}(S)| \geq \exp(-c I(f))$.
\end{conjecture}
\indent It seems that it would be helpful to solve the FMEI conjecture first. \\
\indent Although the two conjectures (if right) help a lot on the understanding of the structure of the Fourier spectrum of the boolean function, the progress on both of them has been slow. (See the introduction for details).\\
\indent We finally give another example on the conjectures as the end of this section.
\begin{example}
\textbf{(And Function)} For any $n>2$, let $f$ be the so-called 'And Function' such that $f(x)=1$ if and only if $x_i=1$ for any $i\in[n]$. It is also written as $x_1 \wedge x_2 ... \wedge x_n$. \\
\indent For such functions, it's not hard to find that $I(f)=\frac{n}{2^{n-1}}$. Also, $\hat{f}(S)=\frac{1}{2^{n-1}}$ for $S \neq \emptyset$ and $\hat{f}(S)=-1+\frac{1}{2^{n-1}}$ for $S = \emptyset$. Thus
\begin{equation}
\begin{aligned}
Ent(f) &= \frac{2^n-1}{2^{2(n-1)}}(2n-2)+(1-\frac{1}{2^{n-1}})^2\log_2 \frac{1}{(1-\frac{1}{2^{n-1}})^2}\\
& \leq \frac{4(n-1)}{2^{n-1}}+2 \log_2 (1+ \frac{1}{2^{n-1}-1})\\
& \leq \frac{4(n-1)}{2^{n-1}}+\frac{2}{2^{n-1}-1} \frac{1}{\log 2}\\
& \leq \frac{4n}{2^{n-1}}.
\end{aligned}
\end{equation}
\indent From the inequality above, we know that $c=4$ is enough for such 'And Function'.
\end{example}
\section{Hypercontractivity and its application}
\indent In this chapter, we'll introduce some classical results in the analysis of boolean functions. In fact, it turns out not easy to study the restrictions behind the Fourier coefficients and the boolean function. Up to now, there are not many useful tools on this problem apart from the hypercontractivity, which we are going to introduce first, and after that, we'll give out some theorems and examples to show its application in the analysis of boolean functions.\\
\indent Actually, one may find out that the ideas behind all the applications of the hypercontractivity in this chapter are similar. We basically use hypercontractivity to bound the low degree Fourier weights and use the 'Markov Argument' to bound the high degree Fourier weights and try to get a contradiction by them. This is perhaps the most powerful and untrivial idea in the analysis of general boolean functions that we already know. 
\subsection{Hypercontractivity and its proof}
\indent We first give the definition of the noise operator.
\begin{definition}
For any $f:\{-1,1\}^n \to \mathbb{R}$, $T_\rho f:=\sum_S \rho^{|S|}\hat{f}(S)X_S$. Here $\rho \in [0,1]$ and $T_\rho$ is called a noise operator.
\end{definition}
\textbf{Remark} One may check that $T_\rho f(x)$ is in fact the expectation of $f(x')$, where $x'$ is given by 'rerandomized' each bit of $x$ with probability $1-\rho$ independently. Such process could be understood as a 'noise', and one may read \cite{5} for further study related to the noise sensitivity of boolean function.\\
\indent The hypercontractivity is stated as follows:
\begin{theorem}
(Hypercontractivity) For any $f:\{-1,1\}^n \to \mathbb{R}$, $\rho \in [0,1]$ we have
\begin{equation}
||T_\rho f||_2 \leq ||f||_{1+\rho^2}.
\end{equation}
\end{theorem}
\textbf{Proof} We first prove the situation for $n=1$. We write $\hat{f}(\emptyset)=a$ and $\hat{f}(\{1\})=b$, so that $f(-1)=a-b$, $f(1)=a+b$, and
\begin{equation}
||T_\rho f||_2=(\rho^2 b^2+ a^2)^{\frac{1}{2}},
\end{equation}
\begin{equation}
||f||_{1+\rho^2}=(\frac{|a-b|^{1+\rho^2}+|a+b|^{1+\rho^2}}{2})^{\frac{1}{1+\rho^2}}.
\end{equation}
\indent Without loss of generality, we may assume that $a\geq 0$ and $b \geq 0$. Since $\rho^2 b^2+ a^2 \leq \rho^2 a^2+ b^2$ if $a < b$, we can also assume that $a \geq b$. At last, it is trivial for the situation of $a=0$, so that for $a\neq 0$, if we write $t=\frac{b}{a}$, then $t\in [0,1]$ and the inequality is equivalent to
\begin{equation}
(1+\rho^2 t^2)^{\frac{1}{2}} \leq (\frac{(1-t)^{1+\rho^2}+(1+t)^{1+\rho^2}}{2})^{\frac{1}{1+\rho^2}}.
\end{equation}
\indent Note that by a binomial expansion, we have
\begin{equation}
\frac{(1-t)^{1+\rho^2}+(1+t)^{1+\rho^2}}{2} \geq 1+\frac{(1+\rho^2)\rho^2}{2}t^2.
\end{equation}
\indent Also, since $\frac{1+\rho^2}{2} \in [0,1]$, we have
\begin{equation}
(1+\rho^2 t^2)^{\frac{1+\rho^2}{2}} \leq 1+\frac{(1+\rho^2)\rho^2}{2}t^2.
\end{equation}
\indent Combine the last two inequalities together we already prove the situation for $n=1$.\\ \\
\indent If the theorem works for $n=m$, we now try to prove it for $n=m+1$. We define $f_1(x)=f(-1,x)$ and $f_2(x)=f(1,x)$ for any $x\in \{-1,1\}^m$. Check the Fourier-Walsh expansion of $f_1+f_2$ and $f_1-f_2$, one may find that
\begin{equation}
||T_\rho f||_2^2=\rho^2||T_\rho (\frac{f_1-f_2}{2})||_2^2+||T_\rho (\frac{f_1+f_2}{2})||_2^2.
\end{equation}
\indent Thus,
\begin{equation}
\begin{aligned}
||T_\rho f||_2^2 & = \frac{1+\rho^2}{4}(||T_\rho f_1||_2^2+||T_\rho f_2||_2^2)+\frac{1-\rho^2}{2}\mathbb{E}_x T_\rho f_1(x) T_\rho f_2(x) \\
& \leq  \frac{1+\rho^2}{4}(||f_1||_{1+\rho^2}^2+||f_2||_{1+\rho^2}^2)+\frac{1-\rho^2}{2}||f_1||_{1+\rho^2}||f_2||_{1+\rho^2}\\
&=\rho^2(\frac{||f_1||_{1+\rho^2}-||f_2||_{1+\rho^2}}{2})^2+(\frac{||f_1||_{1+\rho^2}+||f_2||_{1+\rho^2}}{2})^2\\
& \leq (\frac{1}{2}(||f_1||_{1+\rho^2}^{1+\rho^2}+||f_2||_{1+\rho^2}^{1+\rho^2}))^{\frac{2}{1+\rho^2}}\\
&=||f||_{1+\rho^2}^2.
\end{aligned}
\end{equation}
\indent The first inequality is by the situation for $n=m$ and Cauchy-Schwarz inequality, the second inequality is by the situation for $n=1$. From above we have already proved the situation for $n=m+1$ and furthermore, the whole theorem. \qed \\ \\
\indent \textbf{Remark} The theorem in fact shows that such a transformation between the two normed spaces is a contraction. For a more essential proof, one may read \cite{5}, where the Minkowski inequality is used to show that the tensorization of two such contractions is also a contraction.\\
\subsection{An edge-isoperimetric inequality}
\indent For any boolean function $f$, if we define $\partial_k f(x):=\frac{1}{2}(f(x)-f(\mu_k(x))$, in the preliminaries (the proof of Theorem 2.2.1) we already know that $\partial_k f(x)=\sum_{k\in S} \hat{f}(S)X_S$. On the other hand, $\partial_k f$ only takes value from $\{-1,0,1\}$, so that $||\partial_k f||_{1+\rho^2}^{1+\rho^2}=\mathbb{E}|\partial_k f|=\mathbb{P}(f(x)\neq f(\mu_k(x)))=I_k(f)$, which shows us how the hypercontractivity on $\partial_k f$ help us on the analysis of boolean functions.\\
\indent We now give an example of such technique, the edge-isoperimetric inequality below is also known as KKL Edge-Isoperimetric Theorem sometimes.
\begin{theorem}
There exists constants $c>0$, such that for any boolean function $f$ with $\mathbb{E}f=0$, we have
\begin{equation}
\max_k I_k(f) \geq \exp(-c I(f)).
\end{equation}
\end{theorem}
\textbf{Proof} Note that $||T_\rho \partial_k f||_2^2=\sum_{k\in S} \rho^{2|S|}\hat{f}(S)^2$. Thus by hypercontractivity we have for any $\rho \in [0,1]$
\begin{equation}
\sum_{k\in S}\rho^{2|S|}\hat{f}(S)^2 \leq I_k(f)^{\frac{2}{1+\rho^2}},
\end{equation}
so that
\begin{equation}
\begin{aligned}
\sum_{k\in[n]}I_k(f)^{\frac{2}{1+\rho^2}} & \geq \sum_{S}\rho^{2|S|}|S|\hat{f}(S)^2 \\
&\geq \sum_{|S|\leq 2I(f)}\rho^{2|S|}|S|\hat{f}(S)^2\\
&\geq \sum_{|S|\leq 2I(f)}\rho^{4I(f)}\hat{f}(S)^2.
\end{aligned}
\end{equation}
\indent In the last inequality we use the fact that $\hat{f}(\emptyset)=\mathbb{E}f =0$. On the other hand, we have (sometimes known as the Markov argument)
\begin{equation}
\begin{aligned}
\sum_{|S|>2I(f)}\hat{f}(S)^2 & \leq \frac{1}{2I(f)} \sum_{|S|>2I(f)}|S|\hat{f}(S)^2\\
&\leq \frac{1}{2I(f)} \sum_{S}|S|\hat{f}(S)^2\\
&=\frac{I(f)}{2I(f)}\\
&=\frac{1}{2},
\end{aligned}
\end{equation}
so that $\sum_{|S|\leq 2I(f)}\hat{f}(S)^2 \geq 1-\sum_{|S|>2I(f)}\hat{f}(S)^2 =\frac{1}{2}$. Combine it with (3.17) and take $\rho=\frac{1}{2}$ so that $\frac{2}{1+\rho^2}=\frac{8}{5}$, we have $\max_k \{I_k(f)\}^{\frac{3}{5}}I(f)\geq \sum_{k\in[n]}I_k(f)^{\frac{8}{5}}\geq \frac{1}{2} \cdot 2^{-4I(f)}$. This already shows the inequality we want, as $I(f)=\sum_{S}|S|\hat{f}(S)^2 \geq \sum_{|S|\geq 1}\hat{f}(S)^2=1$.\qed\\ \\
\indent \textbf{Remark} The example that we give in the next section somehow also implies the sharpness of this inequality. When $n$ is large, such an lower bound will be much stronger than the trivial bound $\frac{1}{n}I(f)$, which gives some unexpected restrictions for the value of $I_k(f)$ given $I(f)$. The name 'Isoperimetric' comes from the fact that the influence of a boolean function can also be understood as the perimeter of all the 'pivotal edges' in a hypercube.\\
\subsection{KKL Theorem and tribes function}
\indent In this section, we'll introduce the famous Kahn–Kalai–Linial (KKL) Theorem for boolean functions.
\begin{theorem}
(KKL) There exists a constant $c>0$, such that for any boolean function $f:\{-1,1\}^n \to \{-1,1\}$, there exists $k\in [n]$ such that
\begin{equation}
I_k(f)\geq c Var (f)\cdot \frac{\log n}{n}.
\end{equation}
\end{theorem}
\textbf{Proof} Otherwise $I_k(f)\leq c Var (f)\cdot \frac{\log n}{n}$ for any $k$. Recall that by the hypercontractivity of $\partial_k f$ we have for any $\rho \in [0,1]$, $\sum_{k\in S}\rho^{2|S|}\hat{f}(S)^2 \leq I_k(f)^{\frac{2}{1+\rho^2}}$. Take $\rho=\frac{1}{2}$ we have
\begin{equation}
\begin{aligned}
\sum_{1\leq |S|\leq \frac{\log n}{10}}\hat{f}(S)^2 &\leq \rho^{-2 \frac{\log n}{10}}\sum_{k}I_k(f)^{\frac{2}{1+\rho^2}}\\
&=2^{\frac{1}{5}\log n}\sum_{k}I_k(f)^{\frac{8}{5}}\\
&\leq n^{\frac{1}{5}}\cdot n(c Var(f)\cdot \frac{\log n}{n})^{\frac{8}{5}}\\
&\leq c^{\frac{8}{5}}Var(f)\cdot \frac{(\log n)^{\frac{8}{5}}}{n^{\frac{2}{5}}}.
\end{aligned}
\end{equation}
\indent On the other hand,
\begin{equation}
\sum_{|S|> \frac{\log n}{10}}\hat{f}(S)^2 \leq \frac{10}{\log n} I(f) \leq 10 c Var(f),
\end{equation}
here we got the last inequality by writing $I(f)=\sum_k I_k(f)$. Add the two inequalities together we have
\begin{equation}
\begin{aligned}
Var(f)& =\sum_{|S|\geq 1}\hat{f}(S)^2\\
&\leq (c^{\frac{8}{5}}\cdot \frac{(\log n)^{\frac{8}{5}}}{n^{\frac{2}{5}}}+10 c)Var(f).
\end{aligned}
\end{equation}
Note that $(\log n)^4=o(n)$, so that if we take $c>0$ such that $c^{\frac{8}{5}}\cdot \frac{(\log n)^{\frac{8}{5}}}{n^{\frac{2}{5}}}+10 c<1$ for any $n$, we already have a contradiction. Thus the original theorem is proved. \qed \\\\
\indent \textbf{Remark} The technique in the proof is almost the same as the one in the last section. One may also check that the FMEI conjecture is correct for monotone boolean functions by such technique. We can see that this result is somehow sharp from the example below.\\
\begin{example}
\textbf{(Tribes)} Recall the definition of 'And function'($\wedge$) in Example 1. Similarly we have the 'Or function'($\vee$). The tribes, which are actually boolean functions of dimension $n=Nm$, is defined as $f(x):=\vee_{j=0}^{N-1}(\wedge_{i=1}^{m}x_{jm+i})$, where $N>0$ and $m>0$ is given. Fixed $m>0$ and we take $N=\lfloor \log_{1-2^{-m}}^{\frac{1}{2}}\rfloor$ such that $\mathbb{P}(f(x)=-1)=(1-2^{-m})^N$ is closed to $\frac{1}{2}$. In fact
\begin{equation}
|\mathbb{E}f|\leq 2^{1-m},
\end{equation}
and so that
\begin{equation}
Var(f)\geq 1- 2^{2-2m}.
\end{equation}
Take $m\to \infty$, note that $N=O(-\log (1-2^{-m}))=O(2^m)$, so that for any $k$
\begin{equation}
I_k(f)=(1-2^{-m})^{N-1}\cdot \frac{1}{2^{m-1}}=O(\frac{1}{2^{m-1}})=O(2^{-m}).
\end{equation}
Meanwhile, $Var(f)\cdot \frac{n}{\log n}=O(\frac{\log (Nm)}{Nm})=O(\frac{m}{m\cdot 2^m})=O(2^{-m})$, so that we have for such tribes, $I_k(f)=O(Var(f)\cdot \frac{n}{\log n})$ for any $k$.\\
\end{example}
\indent \textbf{Remark} From the arguments above we can see that tribes are a suitable example for KKL theorem to show it's somehow sharp. In fact furthermore we also have $I(f)=nI_k(f)=O(m)$, so that $I_k(f)=\exp(O(I(f)))$. On the other hand $\mathbb{E}f \to 0$ when $m\to \infty$, which means that such tribes are also 'closed to' an example for Theorem 3.2.1 in the last section.
\subsection{Friedgut's junta theorem}
\indent At the end of this chapter,  we'll give out a beautiful structure theorem for boolean functions. This theorem could be seen as a generalization of KKL Edge-Isoperimetric Theorem and is especially useful when $n>\exp(O(I(f)))$. We'll need the definition of $\epsilon$-concentration first.
\begin{definition}
For $\epsilon>0$, the Fourier spectrum of a boolean function $f$ is $\epsilon$-concentrated on $\mathscr{A}\subseteq \mathscr{F}_n=\{S : S\subseteq [n]\}$ if and only if $\mathbb{P}(\mathcal{S}_f\notin \mathscr{A}) \leq \epsilon$.
\end{definition}
In fact the condition is equivalent to $\sum_{S\notin \mathscr{A}}\hat{f}(S)^2\leq \epsilon$.
\begin{theorem}
(Friedgut's junta theorem) For any boolean function $f$, let $J:=\{k \in [n] : I_k(f)\geq (\frac{\epsilon}{I(f)} \cdot 4^{-\frac{I(f)}{\epsilon}})^{\frac{5}{3}}\}$, we have that the Fourier spectrum of $f$ is $2\epsilon$-concentrated on $\mathscr{A}:=\{S : S \subseteq J , |S| \leq \frac{I(f)}{\epsilon}\}$.
\end{theorem}
\textbf{Proof} By the hypercontractivity of $\partial_k f$ and take $\rho=\frac{1}{2}$ we have $\sum_{k\in S}4^{-|S|}\hat{f}(S)^2\leq I_k^{\frac{8}{5}}$, add all the inequalities for $k\notin J$ we have
\begin{equation}
\begin{aligned}
\sum_{S\cap J \neq \emptyset, |S|\leq\frac{I(f)}{\epsilon}}\hat{f}(S)^2 &\leq 4^{\frac{I(f)}{\epsilon}} \sum_{k\notin J} I_k(f)^{\frac{8}{5}}\\
&\leq 4^{\frac{I(f)}{\epsilon}}\cdot (\max_{k\notin J}I_k(f))^{\frac{3}{5}} \sum_{k\notin J} I_k(f)\\
&\leq \frac{\epsilon}{I(f)}\cdot I(f)\\
&=\epsilon.
\end{aligned}
\end{equation}
For the third inequality we use the definition of $J$ and $\sum_{k\notin J}I_k(f)\leq I(f)$. On the other hand we have $\sum_{|S|>\frac{I(f)}{\epsilon}}\hat{f}(S)^2 \leq \frac{\epsilon}{I(f)} (\sum_S |S| \hat{f}(S)^2)=\epsilon$. Thus the Fourier spectrum of $f$ is $2\epsilon$-concentrated on $\{S : S \subseteq J , |S| \leq \frac{I(f)}{\epsilon}\}$. \qed \\\\
\indent \textbf{Remark} Note that $|J|\leq I(f) (\frac{I(f)}{\epsilon} \cdot 4^{\frac{I(f)}{\epsilon}})^{\frac{5}{3}}= \exp(O(I(f)))$. So for $n>\exp(O(I(f)))$, the Friedgut's junta theorem tells us that there should be lots of bits with little Fourier weights. In fact, we can furthermore show that there exists a boolean function $g$ which only depends on all those $x_k$ with $k\in J$ (this kind of boolean functions are called $|J|$-juntas) such that $||f-g||_2^2 \leq 2\epsilon$. One may read \cite{2} for more information about this.\\\\
\indent With the Friedgut's junta theorem above, we can quickly get the corollary below. 
\begin{corollary}
Given $\delta>0$, there exists a constant $c>0$ such that for any boolean function $f$ with $I(f)>\delta$, we have $\max_S |\hat{f}(S)| \geq \exp(-cI(f)^2)$.
\end{corollary}
\textbf{Proof} Take $\epsilon=\frac{1}{4}$ in the Friedgut's junta theorem. Note that $|J|\leq \exp(O(I(f)))$, thus $|\mathscr{A}|\leq (|J|+1)^{\frac{I(f)}{\epsilon}}\leq \exp(O(I(f)^2))$. On the other hand $\sum_{S\in \mathscr{A}}\hat{f}(S)^2\geq 1-2\epsilon =\frac{1}{2}$, so that $\max_S \hat{f}(S)^2 \geq \frac{1}{2|\mathscr{A}|} \geq \exp(-O(I(f)^2))$, which already implies the result we want. \qed \\\\
\indent \textbf{Remark} One may compare this corollary with the FMEI conjecture. Interested readers can also think about what happens for $0\leq I(f) \leq \delta$ where $\delta$ is small enough. In this case, does this corollary still hold? What about the FMEI conjecture? Does it hold for $0\leq I(f) \leq \delta$?
\section{High-order Influence}
\indent In this chapter, we are going to introduce a new object, named high-order influence, in the analysis of boolean functions. Although it was probably already known for a long time (like in \cite{6}), there are not many results about it in the literature. We'll first give the definitions and then try to develop the theories about it by some applications of it that we've found to show its potentiality.
\subsection{Definitions}
\indent Recall that we have already defined the derivatives of a boolean function $f$ by $\partial_k f(x):=\frac{1}{2}(f(x)-f(\mu_k(x))$ and as a result $I_k(f)=\mathbb{E}(\partial_k f)^2=||\partial_k f||_2^2$. This is, in fact, the situation for the first order and we now give a generalization of it.
\begin{definition}
For $f:\{-1,1\}^n \to \mathbb{R}$, the derivative of it on $S\subseteq [n]$ is defined as $\partial_S f(x):=\frac{1}{2^{|S|}}\sum_{T\subseteq S} (-1)^{|T|}f(\mu_T(x))$. Here $\mu_T(x) \in \{-1,1\}^n$ is defined by flipping the coordinates in $T$ of $x$ such that $\mu_T(x)_i= (-1)^{1_{i\in T}}x_i$. In particular, $\partial_\emptyset f=f$.
\end{definition} 

\begin{definition}
For $f:\{-1,1\}^n \to \mathbb{R}$, the high-order influence (or influence for $|S|=1$) of it on a nonempty set $S\subseteq [n]$ is defined as $I_S(f):=\mathbb{E}(\partial_S f)^2=||\partial_S f||_2^2$.
\end{definition}
\indent We also give the following generalization of the total influence.
\begin{definition}
For $f:\{-1,1\}^n \to \mathbb{R}$ and $1\leq m \leq n$, the total influence of $f$ for the $m^{th}$ degree, is defined as $I_{d=m}(f):=\sum_{|S|=m} I_S(f)$. The total influence of $f$ for the first $m$ degrees, is defined as $I_{d\leq m}(f):=\sum_{1\leq |S|\leq m} I_S(f)$.
\end{definition}
\textbf{Remark} Here the phrases 'order' and 'degree' in some sense refer to the same thing, as $|S|$ can denote both the order of the derivative and the degree of polynomial $X_S$.\\ \\
\indent The following theorem on the high-order influence is well-known.
\begin{theorem}
For any $f:\{-1,1\}^n \to \mathbb{R}$ and $S\subseteq [n]$, $|S| \geq 1$ we have
\begin{equation}
\partial_S f=\sum_{S\subseteq T} \hat{f}(T)X_T
\end{equation}
and
\begin{equation}
I_S(f)=\sum_{S\subseteq T} \hat{f}(T)^2.
\end{equation}
\end{theorem}
\textbf{Proof} We prove the theorem by an induction for $|S|$. If $|S|=1$ and $S=\{k\}$, the equations are $\partial_{\{k\}} f=\sum_{k\in T} \hat{f}(T)X_T$ and $I_{\{k\}}(f)=\sum_{k\in T} \hat{f}(T)^2$, which we already know can be proved by writing the Fourier-Walsh expansion of $f(x)$ and $f(\mu_k(x))$. \\
\indent On the other hand, if the theorem works for $S \setminus \{k\}$, $k \in S$, notice that $\partial_S f=\partial_k(\partial_{S\setminus \{k\}}f)$, since
\begin{equation}
\partial_{S\setminus \{k\}} f=\sum_{S\setminus \{k\}\subseteq T} \hat{f}(T)X_T,
\end{equation}
from the situation of $|S|=1$ we have
\begin{equation}
\partial_k(\partial_{S\setminus \{k\}}f)=\sum_{S\setminus \{k\}\subseteq T, k\in T} \hat{f}(T)X_T.
\end{equation}
Thus we already have the first equation. For the second equation, note that
\begin{equation}
I_S(f)=\mathbb{E}_x(\sum_{S\subseteq T} \hat{f}(T)X_T)^2=\sum_{S\subseteq T} \hat{f}(T)^2,
\end{equation}
which helps us finish the induction.
 \qed \\\\
\indent \textbf{Remark} It is also straightforward that the first equation works for $S=\emptyset$ as well.
\subsection{Second order influence}
\indent Different from other orders, the second order influence of boolean functions has some unique properties of its own. In this section, we'll mostly concentrate on those unique properties.\\
\indent The following two theorems are probably already known as we'll explain in the remark of them below.
\begin{theorem}
For a boolean function $f:\{-1,1\}^n \to \{-1,1\}$, if $k,l\in [n]$, $k \neq l$, we have
\begin{equation}
I_{\{k,l\}}(f)=\mathbb{P}(k,l \in \mathcal{P}(x)).
\end{equation}
Recall that $\mathcal{P}(x)$ is the pivotal set of $x$ for $f$.
\end{theorem}
\textbf{Proof} Notice that
\begin{equation}
\begin{aligned}
I_{\{k,l\}}(f) &=\sum_{k,l\in S} \hat{f}(S)^2\\
&=\mathbb{E}_x(\sum_{k\in S}\hat{f}(S)X_S)(\sum_{l\in S}\hat{f}(S)X_S)\\
&=\mathbb{E}_x\partial_k f(x)\partial_l f(x).
\end{aligned}
\end{equation}
\indent On the other hand, since $f$ is a boolean function, $\partial_k f(x)\partial_l f(x)$ only takes value from $\{0,1\}$,  and $\partial_k f(x)\partial_l f(x)=1$ if and only if $f(x)=-f(\mu_k(x))=-f(\mu_l(x))$. Thus we have
\begin{equation}
\begin{aligned}
I_{\{k,l\}}(f) &=\mathbb{P}(\partial_k f(x)\partial_l f(x)=1)\\
&=\mathbb{P}(f(x)=-f(\mu_k(x))=-f(\mu_l(x)))\\
&=\mathbb{P}(k,l \in \mathcal{P}(x)). 
\end{aligned}
\end{equation}
\qed
\begin{theorem}
For any boolean function $f$, we have the following equations:
\begin{equation}
I(f)=\mathbb{E}_x |\mathcal{P}(x)|=\mathbb{E}|\mathcal{S}_f|.
\end{equation}
\begin{equation}
I_{d=2}(f)=\mathbb{E}_x (\frac{1}{2}|\mathcal{P}(x)|(|\mathcal{P}(x)|-1))=\mathbb{E} (\frac{1}{2}|\mathcal{S}_f|(|\mathcal{S}_f|-1)).
\end{equation}
\begin{equation}
I_{d\leq 2}(f)=\mathbb{E}_x (\frac{1}{2}|\mathcal{P}(x)|(|\mathcal{P}(x)|+1))=\mathbb{E} (\frac{1}{2}|\mathcal{S}_f|(|\mathcal{S}_f|+1)).
\end{equation}
\begin{equation}
2I_{d=2}(f)+I(f)=\mathbb{E}_x |\mathcal{P}(x)|^2=\mathbb{E}|\mathcal{S}_f|^2.
\end{equation}
\end{theorem}
\textbf{Proof} We already know the first equation from Chapter 2. For the second equation, note that by the last theorem 
\begin{equation}
I_{d=2}(f)=\sum_{|S|=2} I_S(f)=\sum_{|S|=2}\mathbb{P}(S\subseteq \mathcal{P}(x))=\mathbb{E}_x (\frac{1}{2}|\mathcal{P}(x)|(|\mathcal{P}(x)|-1)).
\end{equation}
\indent On the other hand,
\begin{equation}
I_{d=2}(f)=\sum_{|S|=2} I_S(f)=\sum_{|S|=2}\sum_{S\subseteq T}\hat{f}(T)^2=\sum_S \frac{1}{2}|S|(|S|-1)\hat{f}(S)^2=\mathbb{E} (\frac{1}{2}|\mathcal{S}_f|(|\mathcal{S}_f|-1)).
\end{equation}
\indent For the last two equations, we simply need to take the sum of the first two equations.\qed \\\\
\indent \textbf{Remark} The first equation is well known. In the introduction of \cite{7}, Garban etc. pointed out the fact that $\mathbb{E}_x |\mathcal{P}(x)|^2=\mathbb{E}|\mathcal{S}_f|^2$ in the fourth equation has been observed by Gil Kalai (in a personal communication). They also point out that it is often the case that $\mathbb{E}|\mathcal{S}_f|^2$ is of the same order of $(\mathbb{E}|\mathcal{S}_f|)^2$. Generally for higher moments, we don't have such equations between like $\mathbb{E}|\mathcal{S}_f|^3$ and $\mathbb{E}_x |\mathcal{P}(x)|^3$.\\ \\
\indent The following inequality is an interesting result implying the restrictions behind the influence and the second order influence of the boolean functions.
\begin{theorem}
For any boolean function $f:\{-1,1\}^n \to \{-1,1\}$, we have
\begin{equation}
\sum_{k\in[n]} -I_k(f)\log I_k(f)\leq 2 \mathbb{E}|\mathcal{S}_f|^2,
\end{equation}
we may assume that $I_k(f) \log I_k(f)=0$ for $I_k(f)=0$ here.
\end{theorem}
\textbf{Proof} Note that by hypercontractivity of $\partial_k f$ we have for any $\rho \in [0,1]$
\begin{equation}
\sum_{k\in S} \rho^{2|S|}\hat{f}(S)^2\leq I_k(f)^{\frac{2}{1+\rho^2}}.
\end{equation}
\indent In fact when $\rho=1$, the inequality above becomes an equation since $I_k(f)=\sum_{k\in S} \hat{f}(S)^2$. If we write $r(\rho)=\sum_{k\in S} \rho^{2|S|}\hat{f}(S)^2-I_k(f)^{\frac{2}{1+\rho^2}}$, we'll have $r(\rho)=0$ for $\rho=1$ and $r(\rho)\leq 0$ for $\rho \in [0,1]$. Thus $r'(1) \geq 0$, which means that when $\rho=1$,
\begin{equation}
\sum_{k\in S} 2|S| \rho^{2|S|-1}\hat{f}(S)^2\geq I_k(f)^{\frac{2}{1+\rho^2}}\cdot (-\frac{4\rho \log I_k(f)}{(1+\rho^2)^2}),
\end{equation}
so that
\begin{equation}
-I_k(f)\log I_k(f) \leq \sum_{k\in S}2|S| \hat{f}(S)^2.
\end{equation}
Take the sum for all $k\in[n]$ we have
\begin{equation}
\sum_{k\in[n]} -I_k(f)\log I_k(f) \leq \sum_{S}2|S|^2 \hat{f}(S)^2=2 \mathbb{E}|\mathcal{S}_f|^2,
\end{equation}
which is the result that we want.\qed \\\\
\indent \textbf{Remark} This inequality is also somehow sharp. We can see it from the example of tribes.
\begin{example}
\textbf{(Tribes)} Recall the definition of tribes in Example 3.3.1. We now focus on the second order influence of $f$ in that example (where $N=O(2^m)$). In fact by $I_{\{k,l\}}(f)=\mathbb{P}(k,l\in \mathcal{P}(x))$, we have 
\begin{equation}
I_{\{k,l\}}(f)=\left \{
\begin{aligned}
&(1-2^{-m})^{N-1}\cdot \frac{1}{2^m}, &\lceil \frac{k}{m}\rceil=\lceil \frac{l}{m}\rceil, \\
&(1-2^{-m})^{N-2}\cdot \frac{1}{2^{2m-1}}, &\lceil \frac{k}{m}\rceil \neq \lceil \frac{l}{m}\rceil.
\end{aligned}
\right.
\end{equation}
Thus,
\begin{equation}
\begin{aligned}
I_{d=2}(f)&=\tbinom{m}{2}N (1-2^{-m})^{N-1}\cdot \frac{1}{2^m}+\tbinom{N}{2}m^2 (1-2^{-m})^{N-2}\cdot \frac{1}{2^{2m-1}}\\
&=O(m^2).
\end{aligned}
\end{equation}
\indent On the other hand, recall that $I(f)=O(m)$, $I_k(f)=O(2^{-m})$ for any $k$, we have
\begin{equation}
\mathbb{E}|\mathcal{S}_f|^2=2I_{d=2}(f)+I(f)=O(m^2)
\end{equation}
and 
\begin{equation}
\sum_k -I_k(f)\log I_k(f)=Nm \cdot O(2^{-m})O(m)=O(m^2).
\end{equation}
So that for such $f$, $\sum_k -I_k(f)\log I_k(f)=O(\mathbb{E}|\mathcal{S}_f|^2)$.\\
\end{example}
\indent In the end we want to mention that for a monotone boolean function $f:\{-1,1\}^n \to \{-1,1\}$ and any $S\subseteq [n]$ such that $|S|=2$, $|\partial_S f|$ only takes value from $0$ and $\frac{1}{2}$. (One may check this.) 
\subsection{Application: On the FMEI Conjecture for Low-degree Boolean Functions}
\indent In this section, we'll give another application of the high-order influence of boolean functions which somehow solved the FMEI conjecture for low-degree boolean functions.\\
\indent We first give the definition of degree.
\begin{definition}
The degree of a boolean function $f:\{-1,1\}^n \to \{-1,1\}$ is defined as
\begin{equation}
\deg (f):=\max_{S\subseteq [n]} \{|S|:\hat{f}(S)\neq 0\}
\end{equation}
\end{definition}
\indent Note that $\deg(f)$ is in fact the degree of the polynomial form: $f=\sum_{S} \hat{f}(S)X_S$ and from the perspective of spectral sample we have $\deg(f)=\max \{|\mathcal{S}_f|\}$.\\\\
\indent The following theorem and its proof give us some characteristics of the Fourier spectrum of low-degree boolean functions. The idea is not far from Friedgut's junta theorem and for the proof we basically use the hypercontractivity of $\partial_S f$.
\begin{theorem}
There exists a constant $c>0$ such that for any boolean function $f$ we have 
\begin{equation}
\max_S |\hat{f}(S)| \geq \exp(-c \deg (f)).
\end{equation}
\end{theorem}
\textbf{Proof} For convenience we write $d_f=\deg(f)$ for short. First we may assume that $d_f \geq 1$ since the situation for $d_f=0$ is trivial. Note that
\begin{equation}
\begin{aligned}
I_{d\leq d_f} (f) & =\sum_{1\leq |S| \leq d_f}\sum_{S\subseteq T}\hat{f}(T)^2\\
&\leq \sum_{T}\sum_{S\subseteq T} \hat{f}(T)^2\\
&= \sum_T 2^{|T|} \hat{f}(T)^2\\
&\leq \sum_T 2^{d_f} \hat{f}(T)^2\\
&=2^{d_f}.
\end{aligned}
\end{equation}
\indent For any $\rho\in [0,1]$, by the hypercontractivity of $\partial_S f$ we have $||T_\rho \partial_S f||_2 \leq ||\partial_S f||_{1+\rho^2}$. On one hand $||T_\rho \partial_S f||_2^2=\sum_{S\subseteq T} \rho^{2|T|}\hat{f}(T)^2$ by the definition, on the other hand since $2^{|S|}\partial_S f \in \mathbb{Z}$ we have $|\partial_S f|^{1+\rho^2}\leq (2^{|S|})^{1-\rho^2} |\partial_S f|^2$, thus $\mathbb{E}_x |\partial_S f|^{1+\rho^2}\leq 2^{(1-\rho^2)|S|}I_S(f)$. Combine all the arguments above and take $\rho=\frac{1}{2}$ we have for any $S$ with $1\leq |S|\leq d_f$,
\begin{equation}
\begin{aligned}
\hat{f}(S)^2 & \leq 2^{2|S|}\cdot ||T_{\frac{1}{2}} \partial_S f||_2^2 \\
&\leq 2^{2d_f}\cdot ||\partial_S f||_{\frac{5}{4}}^2\\
&= 2^{2d_f}\cdot (\mathbb{E}_x |\partial_S f|^{\frac{5}{4}})^{\frac{8}{5}}\\
&\leq 2^{2d_f}\cdot (2^{\frac{3}{4}|S|}I_S(f))^{\frac{8}{5}}\\
&\leq 2^{2d_f} \cdot 2^{\frac{6}{5}d_f} I_S(f)^{\frac{8}{5}}\\
&= 2^{\frac{16}{5}d_f}I_S(f)^{\frac{8}{5}}.
\end{aligned}
\end{equation}
\indent Let $\mathscr{B}:=\{1\leq |S| \leq d_f : I_S(f) \geq 2^{-9d_f}\}$, by the inequality above we have
\begin{equation}
\begin{aligned}
\sum_{S\notin \mathscr{B}\cup \{\emptyset\} }\hat{f}(S)^2 & \leq 2^{\frac{16}{5}d_f} \cdot \sum_{S\notin \mathscr{B}\cup \{\emptyset\}}I_S(f)^{\frac{8}{5}}\\
&\leq 2^{\frac{16}{5}d_f} \cdot (2^{-9d_f})^{\frac{3}{5}} \sum_{|S| \geq 1}I_S(f) \\
&= 2^{-\frac{11}{5}d_f}\cdot  I_{d\leq d_f} (f) \\
&\leq 2^{-\frac{6}{5}d_f}\\
&\leq \frac{1}{2},\\
\end{aligned}
\end{equation}
which means that $\sum_{S\in\mathscr{B}\cup \{\emptyset\}}{f}(S)^2 \geq \frac{1}{2}$. From the definition of $\mathscr{B}$ we have $|\mathscr{B}| \cdot 2^{-9d_f} \leq I_{d\leq d_f} (f)$, so that $|\mathscr{B}| \leq 2^{10d_f} = \exp(O(d_f))$, thus $\max_{S\in\mathscr{B}\cup \{\emptyset\}}\hat{f}(S)^2 \geq \frac{1}{2} (|\mathscr{B}|+1)^{-1}=\exp(-O(d_f))$. This is already enough for the result that we need.\qed\\\\
\indent \textbf{Remark} Note that $\deg (f)=\max |\mathcal{S}_f| \geq\mathbb{E}|\mathcal{S}_f| =I(f)$, so this theorem is strictly weaker than the FMEI conjecture and they are equivalent when $\deg(f)=O(I(f))$. From the example below we can see it's possible that $\deg(f)=I(f)$.\\
\begin{example}
For any $n \in \mathbb{N}^*$, $k\in [n]$, we define the boolean function $f:\{-1,1\}^n\to\{-1,1\}$ by $f(x) :=\prod_{i=1}^k x_i$. It's not hard to find out that $\deg(f)=I(f)=k$.\\
\end{example}
\indent In next section we'll give more examples of low-degree boolean functions.
\subsection{Examples of low-degree boolean functions}
\indent For a boolean function $f:\{-1,1\}^n \to \{-1,1\}$, it is straightforward that $\max_S |\hat{f}(S)| \geq 2^{-n}$, otherwise $\hat{f}(S)=0$ for any $S$ by the definition of $\hat{f}(S)$, which means that Theorem 4.3.1 is trivial when $\deg(f)=O(n)$. For this reason we're interested in examples where $\deg(f)=o(n)$.\\
\indent A large class of such examples are the so-called juntas, where the value of $f$ is only determined by part of the coordinates (and if the value is determined by $o(n)$ coordinates we'll have $\deg(f)=o(n)$). A natural question is that if $f$ is essentially determined by $n$ coordinates, or more rigorously, for any $k\in[n]$, there exists $x\in \{-1,1\}^n$ such that $f(x)\neq f(\mu_k(x))$, is it possible for $\deg(f)=o(n)$?\\
\indent In this section, we'll give some propositions in order to get such examples. Such idea of constructing boolean functions by composition comes from \cite{8}, where they use this idea to give an explicit boolean function which implies the constant in FEI conjecture should be at least 6.278.
\begin{proposition}
Given $n_1, n_2 \in \mathbb{N}^*$ and two boolean functions $f:\{-1,1\}^{n_1}\to \{-1,1\}$, $g:\{-1,1\}^{n_2}\to \{-1,1\}$, if we define boolean function $f_g:\{-1,1\}^{n_1 n_2}\to\{-1,1\}$ by 
\begin{equation}
f_g(x_1, ..., x_{n_1 n_2}):=f(g(x_1, ..., x_{n_2}), g(x_{n_2+1}, ..., x_{2n_2}), ..., g(x_{(n_1-1)n_2+1}, ..., x_{n_1 n_2})),
\end{equation}
then we have
\begin{equation}
\deg(f_g)\leq \deg(f)\deg(g).
\end{equation}
\end{proposition}
\textbf{Proof} We define boolean functions $g_S:\{-1,1\}^{n_1 n_2}\to\{-1,1\}$ where $S\subseteq [n_1]$ by
\begin{equation}
g_S(x):=\prod_{j\in S}g(x_{(j-1)n_2+1}, ..., x_{jn_2}).
\end{equation}
\indent Note that by $f=\sum_S \hat{f}(S) X_S$ we have $f_g=\sum_S \hat{f}(S)g_S$. On one hand it's not hard to find that $\deg (g_S)=|S|\deg(g)$, on the other hand $\hat{f}(S)=0$ for $|S|>\deg (f)$, so we have $\deg (\hat{f}(S)g_S)\leq \deg(f)\deg(g)$ for any $S$. Thus 
\begin{equation}
\deg(f_g)=\deg(\sum_S \hat{f}(S)g_S)\leq \deg(f)\deg(g).
\end{equation} \qed 
\begin{proposition}
Given $n_0 \in \mathbb{N}^*$ and a function $h:\{-1,1\}^{n_0}\to \{-1,1\}$, if we define a series of boolean functions $f^{(k)}:\{-1,1\}^{n_0^k}\to \{-1,1\}$, $k=1,2,...$ inductively by $f^{(1)}:=h$, and $f^{(k+1)}:=f^{(k)}_h$, we have for any $k$
\begin{equation}
\deg(f^{(k)})\leq \deg(h)^k.
\end{equation}
\end{proposition}
\indent Proposition 4.4.2 comes straightforward from Proposition 4.4.1. By a more careful calculation we can actually prove $\deg(f_g)= \deg(f)\deg(g)$ and $\deg(f^{(k)})= \deg(h)^k$.\\
\indent By Proposition 4.4.2 one may find out that if $\deg(h) < n_0$, we've already got a series of boolean function $f^{(k)}:\{-1,1\}^{n_0^k}\to \{-1,1\}$, such that $\deg(f^{(k)})=o(n_0^k)$. One may also check that if furthermore $h$ is essentially determined by $n_0$ coordinates, for any $k\in \mathbb{N}^*$, $f^{(k)}$ is also essentially by $n_0^k$ coordinates. Such $h$ is not hard to find, we give a simple example below.\\
\begin{example} 
Let $h:\{-1,1\}^3 \to \{-1,1\}$ be a boolean function such that $h(x)=-1$ if and only if $x=(1,-1,-1)$ or $x=(-1,1,1)$. One may check that $h(x)$ is essentially determined by all of its coordinates and $\hat{h}(\{1,2,3\})=0$ so that $\deg(h)<3$.
\end{example}

\section{Restrictions}
\indent In this chapter, we'll concentrate on a typical and interesting idea in dealing with boolean functions. The idea is simple: if we fix some coordinates of $x$, we'll find that $f(x)$ turns out to be a lower dimensional boolean function, which might be easier to deal with (we can use induction for example).\\
\indent Similar ideas of restricting certain coordinates have been used for a long time indeed, like in \cite{9} or even in the proof of hypercontractivity, and in \cite{6} they have a more systematic introduction. We first give the definition about the idea and using it, we'll prove a theorem on estimating the spectrum entropy $Ent(f)$, which might be a step towards the proof of FEI conjecture.
\subsection{Definition}
\indent We give the following definition of restricted boolean functions.
\begin{definition}
Given a boolean function $f:\{-1,1\}^n\to \{-1,1\}$ and $J\subseteq [n]$, $x\in \{-1,1\}^n$, we define the restricted boolean function $f_{J^c \to x}:\{-1,1\}^J \to \{-1,1\}$ such that for any $y\in \{-1,1\}^J$, $f_{J^c \to x}(y)=f(z)$, where $z\in \{-1,1\}^n$ and $z_i=1_{i\in J}\cdot y_i+1_{i\notin J}\cdot x_i$ for any $i\in [n]$.
\end{definition}
\indent As before, we assume that $x$ obeys the uniform distribution on $\{-1,1\}^n$ if we consider it as a random vector. (Under such assumption, $f_{J^c \to x}$ are also called random restrictions in \cite{6}.)\\
\indent With the definition above, we can actually get numerous equations between the value and Fourier coefficients of a boolean function and the restricted boolean function of it. Here we give out one of them, which will be useful in proving the theorems in the next section.
\begin{theorem}
For any boolean function $f:\{-1,1\}^n \to \{-1,1\}$ and $k$, $J$ such that $k\in J \subseteq [n]$, we have
\begin{equation}
\mathbb{E}_x \sum_{k\in S\subseteq J} \widehat{f_{J^c \to x}}(S)^2 =I_k(f).
\end{equation}
\end{theorem}
\textbf{Proof} Note that $\sum_{k \in S\subseteq J} \widehat{f_{J^c \to x}}(S)^2=I_k(f_{J^c \to x})$ and 
\begin{equation}
\begin{aligned}
\mathbb{E}_x I_k(f_{J^c \to x}) &=\mathbb{E}_x \mathbb{E}|\partial_k f_{J^c \to x}|\\
&=\mathbb{E}|\partial_k f|\\
&=I_k(f),
\end{aligned}
\end{equation}
thus the equation we need is straightforward.\qed
\subsection{Estimating the spectrum entropy of boolean functions}
\indent In this section, we'll try to estimate the spectrum entropy $Ent(f)$ using the restricted boolean functions.\\
\indent Given a boolean function $f:\{-1,1\}^n\to \{-1,1\}$, for any $V \subseteq [n]$, $\epsilon \in [0, \frac{1}{2})$, we define the moments of restricted Fourier coefficients as
\begin{equation}
M_{V,\epsilon}(f):=\mathbb{E}_{x \in \{-1,1\}^n}\sum_{S\subseteq V}{|\widehat{f_{V^{c} \to x}}(S)|^{2(1+\epsilon)}}.
\end{equation} 
Note that we have $M_{V,0}(f)=1$ directly from the definition. On the other hand $M_{[n],\epsilon}(f)=\sum_S |\hat{f}(S)|^{2(1+\epsilon)}$, so that $Ent(f)=-\frac{1}{ln2} \frac{dM_{[n],\epsilon}(f)}{d\epsilon}(0)$. We'll try to bound $M_{V,\epsilon}(f)$ and get the bound of $Ent(f)$ by calculating the derivative of $M_{[n],\epsilon}(f)$ with respect to $\epsilon$.\\
\indent Under this motivation, we give the theorem below:
\begin{theorem}
Given a boolean function $f:\{-1,1\}^n \to \{-1,1\}$, for any $V_1 \subset [n]$, $k \in [n]$ but $k \notin V_1$, $\epsilon \in (0, \frac{1}{2})$, if we write $V_2=V_1 \cup \{k\}$, we have
\begin{equation}
M_{V_2,\epsilon}(f)-M_{V_1,\epsilon}(f) \geq -I_k(f)(3\epsilon +2\epsilon^2+ (\frac{I_k(f)}{4})^{-\epsilon}-1).
\end{equation}
\end{theorem}
\indent For the proof of the theorem, we'll need the following lemma:
\begin{lemma}
For $0 \leq a \leq b \leq 1$, $\epsilon \in (0,\frac{1}{2})$, we have
\begin{equation}
\frac{(\sqrt{b}+\sqrt{a})^{2(1+\epsilon)}+(\sqrt{b}-\sqrt{a})^{2(1+\epsilon)}}{2}-a^{1+\epsilon}-b^{1+\epsilon} \leq (3\epsilon +2\epsilon^2)a + (b^\epsilon-a^\epsilon)a.
\end{equation}
\end{lemma}
\textbf{Proof} Note that $\binom{2+2\epsilon}{2m}<0$ for $m\in \mathbb{N}$, $m\geq2$, by a binomial expansion we have
\begin{equation}
\begin{aligned}
 \frac{(\sqrt{b}+\sqrt{a})^{2(1+\epsilon)}+(\sqrt{b}-\sqrt{a})^{2(1+\epsilon)}}{2} &\leq b^{1+\epsilon}+ \binom{2+2\epsilon}{2}b^{\epsilon}a\\
&= b^{1+\epsilon}+ (1+3\epsilon +2\epsilon^2)b^{\epsilon}a\\
&\leq b^{1+\epsilon}+(3\epsilon +2\epsilon^2)a+b^{\epsilon}a.
\end{aligned}
\end{equation}
Thus the inequality we need is straightforward.\qed\\\\
\indent Now we can focus on the proof of Theorem 5.2.1.\\
\indent \textbf{Proof of Theorem 5.2.1} In this proof we sometimes write $I_k$ which means $I_k(f)$ for short. Note that we have $M_{V_1,\epsilon}(f)=\mathbb{E}_x\sum_{S\subseteq V_1}{|\widehat{f_{V_1^{c} \to x}}(S)|^{2(1+\epsilon)}}=\mathbb{E}_x\sum_{S\subseteq V_1}{|\widehat{f_{V_1^{c} \to \mu_k(x)}}(S)|^{2(1+\epsilon)}}$, so that 
\begin{equation}
M_{V_1,\epsilon}(f)=\mathbb{E}_x\frac{1}{2}\sum_{S\subseteq V_1}{(|\widehat{f_{V_1^{c} \to x}}(S)|^{2(1+\epsilon)}+|\widehat{f_{V_1^{c} \to \mu_k(x)}}(S)|^{2(1+\epsilon)}}).
\end{equation}
Thus
\begin{equation}
\begin{aligned}
M_{V_2,\epsilon}(f)-M_{V_1,\epsilon}(f)= & -\mathbb{E}_x\sum_{S\subseteq V_1}{(\frac{1}{2}(|\widehat{f_{V_1^{c} \to x}}(S)|^{2(1+\epsilon)}+|\widehat{f_{V_1^{c} \to \mu_k(x)}}(S)|^{2(1+\epsilon)}}) \\
&-|\widehat{f_{V_2^{c} \to x}}(S)|^{2(1+\epsilon)}-|\widehat{f_{V_2^{c} \to x}}(S\cup \{k\})|^{2(1+\epsilon)}).
\end{aligned}
\end{equation}
\indent On the other hand, it's not hard to find that $|\widehat{f_{V_1^{c} \to x}}(S)|$ and $|\widehat{f_{V_1^{c} \to \mu_k(x)}}(S)|$ takes value from $|\widehat{f_{V_2^{c} \to x}}(S)+\widehat{f_{V_2^{c} \to x}}(S\cup \{k\})|$ and $|\widehat{f_{V_2^{c} \to x}}(S)-\widehat{f_{V_2^{c} \to x}}(S\cup \{k\})|$ separately (the order might be changed) from the definition. If we write $a_{x,S}=\min\{\widehat{f_{V_2^{c} \to x}}(S)^2,\widehat{f_{V_2^{c} \to x}}(S\cup \{k\})^2\}$ and $b_{x,S}=\max\{\widehat{f_{V_2^{c} \to x}}(S)^2,\widehat{f_{V_2^{c} \to x}}(S\cup \{k\})^2\}$, we'll have $0\leq a_{x,S} \leq b_{x,S} \leq 1$ and
\begin{equation}
M_{V_2,\epsilon}(f)-M_{V_1,\epsilon}(f)= -\mathbb{E}_x\sum_{S\subseteq V_1}{(\frac{1}{2}((\sqrt{b_{x,S}}+\sqrt{a_{x,S}})^{2(1+\epsilon)}+(\sqrt{b_{x,S}}-\sqrt{a_{x,S}})^{2(1+\epsilon)})-a_{x,S}^{1+\epsilon}-b_{x,S}^{1+\epsilon}}).
\end{equation}
Here we can use Lemma 5.2.1 and get that
\begin{equation}
M_{V_2,\epsilon}(f)-M_{V_1,\epsilon}(f) \geq -\mathbb{E}_x\sum_{S\subseteq V_1}((3\epsilon+2\epsilon^2)a_{x,S}+(b_{x,S}^\epsilon-a_{x,S}^\epsilon)a_{x,S}).
\end{equation}
Another fact is that by Theorem 5.1.1 we have 
\begin{equation}
\mathbb{E}_x\sum_{S\subseteq V_1}a_{x,S} \leq \mathbb{E}_x\sum_{S\subseteq V_1}\widehat{f_{V_2^{c} \to x}}(S\cup \{k\})^2=I_k,
\end{equation}
so that we only need to bound $\mathbb{E}_x\sum_{S\subseteq V_1}(b_{x,S}^\epsilon-a_{x,S}^\epsilon)a_{x,S}$.\\\\
\indent For this purpose,  we use Hölder's inequality and notice that $\mathbb{E}_x\sum_{S\subseteq V_1} b_{x,S} \leq 1$, $\frac{1}{1-\epsilon} \geq 1+\epsilon$ to get 
\begin{equation}
\begin{aligned}
\mathbb{E}_x\sum_{S\subseteq V_1}b_{x,S}^\epsilon a_{x,S} & \leq (\mathbb{E}_x\sum_{S\subseteq V_1} b_{x,S})^\epsilon (\mathbb{E}_x\sum_{S\subseteq V_1}a_{x,S}^{\frac{1}{1-\epsilon}})^{1-\epsilon}\\
&\leq (\mathbb{E}_x\sum_{S\subseteq V_1}a_{x,S}^{1+\epsilon})^{1-\epsilon}.
\end{aligned}
\end{equation}
\indent We suppose $A=\mathbb{E}_x\sum_{S\subseteq V_1}a_{x,S}^{1+\epsilon}$, thus $\mathbb{E}_x\sum_{S\subseteq V_1}(b_{x,S}^\epsilon-a_{x,S}^\epsilon)a_{x,S}\leq A^{1-\epsilon}-A$. Note that $A \leq \mathbb{E}_x\sum_{S\subseteq V_1}a_{x,S}\leq I_k$, and $\frac{d(A^{1-\epsilon}-A)}{dA}=0$ if and only if $A=(1-\epsilon)^{\frac{1}{\epsilon}}\geq \frac{1}{4} \geq \frac{1}{4}I_k$, so that the maximal point of $A^{1-\epsilon}-A$ is in $[\frac{1}{4}I_k, I_k]$, which means that $\mathbb{E}_x\sum_{S\subseteq V_1}(b_{x,S}^\epsilon-a_{x,S}^\epsilon)a_{x,S}\leq A(A^{-\epsilon} -1) \leq I_k((\frac{I_k}{4})^{-\epsilon}-1)$. Combine it with (5.69) and (5.70) we directly get the inequality we want to prove. \qed \\\\
\indent By Theorem 5.2.1, we can get the following results, which give rather good estimation for $Ent(f)$.
\begin{theorem}
There exist $c_1, c_2>0$, such that for any boolean function $f:\{-1,1\}^n\to \{-1,1\}$ we have
\begin{equation}
Ent(f)\leq c_1 I(f)+ c_2 \sum_{k\in[n]}-I_k(f)\log I_k(f).
\end{equation}
\end{theorem}
\textbf{Proof} Take $V_1=[k-1]$ and $V_2=[k]$ where $k=1,2,...,n$ in Theorem 5.2.1, we have for any $k\in[n]$, $\epsilon \in (0,\frac{1}{2})$,
\begin{equation}
M_{[k],\epsilon}(f)-M_{[k-1],\epsilon}(f) \geq -I_k(f)(3\epsilon +2\epsilon^2+ (\frac{I_k(f)}{4})^{-\epsilon}-1).
\end{equation}
Add all the inequalities for each $k\in [n]$ together, note that $M_{\emptyset,\epsilon}(f)=\mathbb{E}_x |f(x)|^{2(1+\epsilon)}=1$ we have for any $\epsilon \in (0,\frac{1}{2})$,
\begin{equation}
M_{[n],\epsilon}(f) \geq 1-(3\epsilon +2\epsilon^2)I(f)- \sum_{k=1}^n ((\frac{I_k(f)}{4})^{-\epsilon}-1)I_k(f).
\end{equation}
We also have $M_{[n], 0}(f)=1$, thus
\begin{equation}
\begin{aligned}
Ent(f) &= -\frac{1}{\ln2} \frac{dM_{[n],\epsilon}}{d\epsilon}(0) \\
&\leq \frac{1}{\ln2}\frac{d((3\epsilon +2\epsilon^2)I(f)+ \sum_{k=1}^n ((\frac{I_k(f)}{4})^{-\epsilon}-1)I_k(f))}{d\epsilon}(0)\\
&= \frac{1}{\ln2}(3I(f)+\sum_{k=1}^n{I_k(f) \log\frac{4}{I_k(f)}}),
\end{aligned}
\end{equation}
which is the result that we need.\qed\\
\indent \textbf{Remark} As before we may assume that $I_k(f)\log I_k(f)=0$ when $I_k(f)=0$,  since $I_k(f)\in [0,1]$, $-I_k(f)\log I_k(f)$ is always non-negative. On the other hand, there seems to be no obvious relationship between $\sum_k -I_k(f)\log I_k(f)$ and $I(f)$ or $I(f)\log I(f)$, $I(f)^2$, etc.
\begin{theorem}
There exists a constant $c>0$, such that for any boolean function $f:\{-1,1\}^n\to \{-1,1\}$ we have
\begin{equation}
Ent(f)\leq c \cdot \mathbb{E}|\mathcal{S}_f|^2.
\end{equation}
\end{theorem}
\textbf{Proof} Note that $I(f)=\mathbb{E}|\mathcal{S}_f|\leq \mathbb{E}|\mathcal{S}_f|^2$ and by Theorem 4.2.3 we have 
\begin{equation}
\sum_{k\in[n]} -I_k(f)\log I_k(f)\leq 2 \mathbb{E}|\mathcal{S}_f|^2,
\end{equation}
thus we directly get the result that we need from Theorem 5.2.2.\qed\\\\
\indent \textbf{Remark} By Theorem 4.2.2 we also know that $Ent(f)\leq O(I_{d\leq 2}(f))$ and $Ent(f) \leq O(\mathbb{E}_x|\mathcal{P}(x)|^2)$ for any boolean function $f$. We also want to repeat here that it is often the case that $\mathbb{E}|\mathcal{S}_f|^2$ is of the same order of $(\mathbb{E}|\mathcal{S}_f|)^2=I(f)^2$. From Example 2.4.1, where $\mathbb{E}|\mathcal{S}_f|^2=\mathbb{E}|\mathcal{S}_f|=1$ and $Ent(f)=\log_2 n$, one may find that this theorem doesn't always hold for non-boolean functions as well. This theorem is actually also a corollary for Theorem 1.3 in \cite{6}.\\\\
\emph{Acknowledgement: the author would like to thank Prof. Juhan Aru for his encouragements and valuable advices and thank for the beauty of boolean functions. }

\end{document}